\documentclass[preprintnumbers,amsmath,amssymb]{revtex4}
\usepackage{color}
\usepackage{graphicx}
\usepackage{epsfig}
\usepackage{epstopdf,hyperref}
\usepackage{bm}
\newtheorem{corollary}{Corollary}[section]
\newtheorem{remark}{Remark}[section]
\newtheorem{definition}{Definition}[section]
\newtheorem{theorem}{Theorem}[section]
\newtheorem{proposition}{Proposition}[section]
\newtheorem{lemma}{Lemma}[section]

\begin{document}
\title{Dynamics of nonlocal and local discrete Ginzburg-Landau equations: global attractors and their congruence}
	
\author{Dirk Hennig}
	\author{Nikos I. Karachalios}
\affiliation{Department of Mathematics, University of Thessaly, Lamia GR35100, Greece}
\date{\today}

\begin{abstract}
Discrete Ginzburg-Landau (DGL) equations with non-local nonlinearities have been established as significant inherently discrete models in numerous physical contexts, similar to their counterparts with local nonlinear terms.  We study two prototypical examples of non-local and local DGLs on the one-dimensional infinite lattice. For the non-local DGL, we identify distinct scenarios for the asymptotic behavior of the globally existing in time solutions depending on certain parametric regimes. One of these scenarios is associated with a restricted compact attractor according to J. K. Hale's definition.  We also prove the closeness of the solutions of the two models  in the sense of a "continuous dependence on their initial data" in the $l^2$ metric under general conditions on the intrinsic linear gain or loss incorporated in the model. As a consequence of the closeness results, in the dissipative regime we establish the congruence of the attractors  possessed  by the semiflows of the non-local and of the local model respectively, for initial conditions in a suitable domain of attraction defined by the non-local system. 
\end{abstract}

\maketitle
\section{Introduction}

\noindent 
In this paper, we study the following two discrete complex Ginzburg-Landau equations: 
\begin{equation}
\dot{u}_n=(1-\delta)u_n+(1+i\alpha)
(u_{n+1}-2u_n+u_{n-1})-(1+i\beta)\,|u_{n}|^2u_n+g_n,\,\,\,n\in {\mathbb{Z}},\label{eq:DNLS}
\end{equation}
with $u_n \in {\mathbb{C}}$ and initial conditions:
\begin{equation}
 u_{n}(0)=u_{n,0},\,\,\,n \in {\mathbb{Z}},\label{eq:icsDNLS}
\end{equation}
and 
\begin{equation}
\dot{v}_n=(1-\delta)v_n+(1+i\alpha)(v_{n+1}-2v_n+v_{n-1})-\frac{1}{2}(1+i\beta)|v_n|^2(v_{n+1}+v_{n-1}) +g_n,\,\,\,n\in {\mathbb{Z}},\label{eq:ALS}
\end{equation}
with $v_n \in {\mathbb{C}}$ and initial conditions:
\begin{equation}
 v_{n}(0)=v_{n,0},\,\,\,n \in {\mathbb{Z}}.\label{eq:icsALS}
\end{equation}
In what follows we refer to the systems (\ref{eq:DNLS}) and (\ref{eq:ALS}) as the {\it L-DGL} and {\it NL-DGL}, respectively. For real parameter $\delta >0$, the first term on the right-hand side of (\ref{eq:DNLS}) and (\ref{eq:ALS}) is of  intrinsic linear gain [loss (damping)] nature if $\delta<1$ ($\delta>1$), while $g=(g_n)_{n\in {\mathbb{Z}}}\neq 0$ acts in both systems as a general external force. Further, we assume $\alpha,\beta \in {\mathbb{R}}_+$. 

Nonlinear lattices of the form (\ref{eq:DNLS}) or (\ref{eq:ALS}), are encountered in several diverse branches of physics, ranging from the dynamics of coupled laser systems and waveguides in 1D \cite{UGL1},\cite{UGL2},\cite{BUGL2},  and 2D-setups \cite{UGL3},\cite{UGL4},\cite{BUGL4}, to the description of the dynamical behavior of low dimensional fluid dynamical systems \cite{UGL5}. For instance, models of the form (\ref{eq:ALS}) with nonlocal nonlinearities determine the evolution of the amplitude of the separatrices between a linear array of vortices or coupled waves observed in hydrodynamic systems of low dimensionality \cite{AGL1},\cite{AGL2}.  Furthermore, the NL-DGL \eqref{eq:ALS}  and other variants with generalized nonlocal nonlinearities have proven to be relevant models for the construction of soliton solutions in the presence of gain and loss \cite{AGL3},\cite{AGL4} and investigations of their potential stability. In the above physical contexts, other relevant discrete DGL equations involve saturable nonlinearities \cite{UGL6}. 

Notice that in (\ref{eq:DNLS}) and (\ref{eq:ALS}) the nonlinear terms are both of cubic order, and an important common feature of both systems, is that in 
the continuum limit approximate the 
Ginzburg-Landau partial differential equation (PDE),
\begin{equation}
\partial_t u=(1-\delta)u+(1+i\alpha)
\partial_{xx} u -(1+i\beta)\,|u|^2u+g(x),\label{eq:CGL}
\end{equation} 
which is a spatially driven counterpart of the complex Ginzburg-Landau equation; a spatial driver may appear after application of a usual gauge-time transformation in the case of a spatiotemporal driving force with harmonic in time driving  \cite{dGLE1}-\cite{dGLE5}.
However, in the discrete regime, the systems (\ref{eq:DNLS}) and (\ref{eq:ALS})  {\em are markedly different} in the sense that the nonlinear terms in (\ref{eq:ALS}) are of nonlocal nature compared to the local terms in (\ref{eq:DNLS}).

Herein, we compare the solution properties of the two infinite lattice dynamical systems presented by (\ref{eq:DNLS}) and (\ref{eq:ALS}) with particular attention to their global  asymptotic behaviour. 
We introduce the distance variable 
$\Delta =(\Delta_n)_{n\in \mathbb{Z}}$ where $\Delta_n=u_n-v_n$. The main results are the following:
\begin{enumerate}
	\item {\em Closeness of the solutions}.  Solutions of both systems remain close when initiated from sufficiently close initial data according to 
\begin{theorem}
	\label{Theorem:closeness}	
	Assume that the initial condition $(u_n(0))_{n\in {\mathbb{Z}}}=u^0$ of the L-DGL (\ref{eq:DNLS})-(\ref{eq:icsDNLS}) and the initial condition $(v_n(0))_{n\in {\mathbb{Z}}}=v^0$ of the NL-DGL  (\ref{eq:ALS})-(\ref{eq:icsALS}) satisfy 
	\begin{eqnarray}
	|| u(0)-v(0)||_{l^2}&\le& C_0 \epsilon^3,\label{eq:distance0}\\
	|| u(0)||_{l^2}&\le& C_{u,0}\epsilon:=R_u(\epsilon),\,\,\, || v(0)||_{l^2}\le C_{v,0}\epsilon:=R_v(\epsilon),\label{eq:initialuv}
	\end{eqnarray}
	for some positive constants $C_{u,0}$, $C_{v,0}>0$. Then, we have the following cases for the closeness of the corresponding solutions:
	\begin{enumerate}
		\item Case $\delta>1$.  For any $t\in (0, \infty]$, the corresponding solutions to Eqs. (\ref{eq:DNLS}) and (\ref{eq:ALS}) fulfill  the estimate
		\begin{equation}
		|| u(t)-v(t)||_{l^2}\le C \epsilon^3,\label{eq:boundy}
		\end{equation}
		where the constant $C(C_0,C_{u,0},C_{v,0},\beta, \delta)$ is independent of $t$. Furthermore, there exists a constant $C_1(C_0,C_{u,0},C_{v,0},\beta, \delta)$, still independent of $t$, such that
		\begin{eqnarray}
		\label{cd1}
		\limsup_{t\rightarrow\infty}||\Delta(t)||_{l^2}\leq C_1\epsilon^3.
		\end{eqnarray}
		\item Case $\delta\leq 1$.  For any arbitrary $0<T_f<\infty$ there exists a constant $C_2(C_0,C_{u,0},C_{v,0},\beta, \delta,T_f)$ such that the the corresponding solutions to Eqs. (\ref{eq:DNLS}) and (\ref{eq:ALS}) fulfill  the estimate
		\begin{equation}
		\label{eq:boundylin}
		|| u(t)-v(t)||_{l^2}\le C_2\epsilon^3.
		\end{equation}
	\end{enumerate}
\end{theorem}	
A first application of Theorem \ref{Theorem:closeness} is in the context of {\em dissipative discrete solitons}, mentioned above:  Localised solutions have been constructed for the NL-DGL \eqref{eq:ALS},\cite{AGL3},\cite{AGL4},\cite{Abdul2} as deformed continuations of the exact soliton solutions of the integrable Ablowitz-Ladik lattice \cite{AL}, \cite{Ablowitz}. Theorem \ref{Theorem:closeness} rigorously establishes that the L-DGL model \eqref{eq:DNLS} admits $\mathcal{O}(\epsilon)$ solutions that stay $\mathcal{O}(\epsilon^3)$ close to the solution of the NL-DGL \eqref{eq:ALS}.  In the undamped case $\delta=1$ we prove that the growth of the distance $||\Delta(t)||_{l^2}$ is at most linear for $t\in [0, T_f]$, while in the case of the intrinsic forcing $\delta<1$ we found that the growth is at most exponential; the latter case is similar to the exponential growth of the distance function when the inviscid limit of the complex Ginzburg-Landau PDE the Nonlinear Schr\"odinger PDE is considered \cite{JWU},\cite{OG}.
\item {\em Congruence of attractors in the dissipative regime $\delta>1$}.  A second application of Theorem \ref{Theorem:closeness} is in the context of infinite dimensional dynamical systems \cite{Hale},\cite{Temam},\cite{Robinson},\cite{Chueshov}  concerning the existence of attractors possessed by the L-DGL and NL-DGL  systems. 

It is useful for completeness and in order to highlight their differences,  to recall the definitions of the global and restricted attractor from  \cite{Hale}, stated for a generic dynamical system in $l^2$ denoted by $\left\{S(t)\right\}_{t\ge 0}$ representing one of the systems \eqref{eq:DNLS} and \eqref{eq:ALS}.
%
\begin{definition} 
	\label{GAu1}	
	A set ${\cal A}\subset l^2$ is called a global attractor for the semigroup  $\left\{S(t)\right\}_{t\ge 0}$  in $l^2$ if 
	(i)  ${\cal A}\neq 0$ is a compact subset of $l^2$, (ii) an invariant set, that is,  
	$S(t){\cal A}={\cal A}$ for all $t\ge 0$, 
	and, (iii) an attracting set for $\left\{S(t)\right\}_{t\ge 0}$ in $l^2$, viz. for all bounded $\mathcal{B}\subset l^2$, it holds that 
	$\lim_{t\rightarrow \infty}\,dist(S(t)\mathcal{B},{\cal A})=0$, where the Hausdorff semi-distance between two nonempty subsets $\mathcal{U},\mathcal{V}$ of $l^2$ is determined by
	\begin{equation}
	dist(\mathcal{U},\mathcal{V})=\sup_{u\in \mathcal{U}}\,\inf_{v\in \mathcal{V}}\,d(u,v)_{l^2}.
	\end{equation}
\end{definition}

\begin{definition}
	\label{RGAu3}
We say that ${\cal A}_u\subset l^2$ is a restricted global attractor for  $\left\{S(t)\right\}_{t\ge 0}$
	in $l^2$ if for some closed, nonempty subset $\mathcal{U}$
	of $l^2$, $S(t):\mathcal{U}\rightarrow \mathcal{U}$ ($t\ge 0$) is a  semigroup on $\mathcal{U}$ s.t. ${\cal A}_u$ is a global attractor for $\left\{S(t)\right\}_{t\ge 0}$
	restricted to $\mathcal{U}$, viz. (i) $S(t){\cal A}_u={\cal A}_u$ for $t\ge 0$, (ii) ${\cal A}_u$ is compact, and (iii) ${\cal A}_u$ attracts solutions of bounded subsets of $\mathcal{U}$.
\end{definition}

When $\delta>1$ the constant $C$ in the closeness of solutions estimate \eqref{eq:boundy} is independent of $t$ and moreover, as $t\rightarrow\infty$, the distance function becomes uniformly bounded by $\mathcal{O}(\epsilon^3)$, as dictated by the limit \eqref{cd1}.  This is the key element in proving that the global attractor $\mathcal{A}_u$ belonging to the dynamical system L-DGL and the {\em restricted attractor $\mathcal{A}_v$} in the sense of \cite{Hale} belonging to the dynamical system NL-DGL coincide for a suitable choice of their initial conditions: Let $\mathcal{B}_{v,0}$  denote the absorbing set for restricted initial data $v^0$ (i.e. satisfying a certain smallness condition) of the dynamical system  NL-DGL. We prove
\begin{theorem}
	Let $\delta>1$ and assume that the initial conditions $u^0$ of the L-DGL (\ref{eq:DNLS})-(\ref{eq:icsDNLS}) and the initial condition $v^0$ of the NL-DGL  (\ref{eq:ALS})-(\ref{eq:icsALS}) satisfy  \eqref{eq:initialuv}  and are such that $u^0,\,v^0\in \mathcal{B}_{v,0}$ . Then the  attractors
	${\cal A}_u$ and ${\cal A}_v$ in $\ell^2$ coincide according to
	\begin{equation}
	\label{distS}
	{\rm dist}\left({\cal A}_v,{\cal A}_u\right)=0. 
	\end{equation}\label{Theorem:congruence}
\end{theorem}
\end{enumerate}
In \eqref{distS}, the distance between the two sets is the usual Haussdorff semi-distance \cite{Temam}. 

On route to the proof of Theorem \ref{Theorem:congruence} we shall discover much richer scenarios for the dynamics of the NL-DGL concerning the global existence of solutions and invariant sets than for the L-DGL, highlighting the impact of the nonlocal nonlinearity. For instance, a careful manipulation of a Ricatti differential inequality (instead of  Gronwall's lemma and its variants) reveals the emergence of distinct regimes for the existence of global in time solutions for all initial data, the existence of an invariant annulus with a potential convergence to a particular state, the existence of a non-escaping ball depending on the size of the initial data, as well as the existence of the restricted absorbing ball giving rise to the restricted attractor. To establish the compactness of the attracting sets we employ an alternative contradiction argument to prove the asymptotic tail property instead of the tail-estimates method \cite{Bates},\cite{Zhou0},\cite{Zhou},\cite{Zhou1},\cite{Zhou3},\cite{Wang},\cite{Zhou2},\cite{Nikos},\cite{NAE2007},\cite{Du}; the tail-estimates method has been established as  a fundamental, generic approach to prove asymptotic compactness not only in autonomous but also in non-autonomous, delayed  or stochastic lattice dynamical systems, see \cite{Stoch1},\cite{Stoch2},\cite{Stoch3},\cite{Stoch4} and references therein. In the case of the NL-DGL \eqref{eq:ALS}, this contradiction argument is combined with generalised comparison principles for ordinary differential equations \cite{McNabb1986}.
 
The presentation of the paper is as follows:  Section  \ref{subsection:existence} contains global existence results and asymptotic estimates for the solutions for the considered GLEs. For local existence of solutions, we combine both systems in a physically significant dissipative version of the Salerno model \cite{Salerno} or Integrable-Nonintegrable-DNLS (IN-DNLS) \cite{Cai94}. For global existence each of the L-DGL and NL-DGL systems is treated separately, revealing in the case of the NL-DGL the relevant distinct parametric regimes and the associated asymptotic estimates.  Section \ref {section:attractors} contains the results on the existence of attracting sets, while section \ref{section:closeness} contains the proofs of the main results, Theorems \ref{Theorem:closeness} and \ref{Theorem:congruence}. Section \ref{SecCon} summarizes the findings and provide a brief  plan for further relevant studies still to be  performed.
\section{Global existence and asymptotic estimates for the solutions}\label{subsection:existence}
\setcounter{equation}{0}
For the current study of existence and uniqueness of a global solution of the NL-DGL and L-DGL, we combine them in a single system
\begin{equation}
\dot{u}_n-(1+i\alpha)(u_{n+1}+u_{n-1}-2u_n)=(1-\delta)u_n-
(1+i\beta)\left[\gamma|u_n|^2u_n+\frac{\mu}{2}(u_{n+1}+u_{n-1})|u_n|^2\right]+g_n,\,\,\,n\in {\mathbb{Z}},\label{eq:systemglobal}
\end{equation}
with $u_n \in {\mathbb{C}}$ and initial conditions:
\begin{equation}
u_{n}(0)=u_{n,0},\,\,\,n \in {\mathbb{Z}}.\label{eq:icsglobal}
\end{equation}
Note that for $\gamma=0$ and $\mu=1$ ($\mu=0$ and $\gamma=1$) the L-DGL (NL-DGL) results from (\ref{eq:systemglobal}). We remark that the combined system \eqref{eq:systemglobal}
is of physical significance itself as it can be considered as an extended dissipative and driven version of the Salerno model \cite{Salerno} or Integrable-Nonintegrable-DNLS (IN-DNLS) \cite{Cai94}. In fact, the latter is obtained from (\ref{eq:systemglobal}) in the limit $\alpha,\beta \rightarrow \infty$, $\delta=0$ and $g=0$. 
\paragraph{Preliminaries.}
\label{parpre}
The functional setting of the initial value problem for the infinite system of ordinary differential equations in (\ref{eq:systemglobal})-(\ref{eq:icsglobal}) consists of the sequence spaces 
\begin{equation}
l^p=\left\{ u=(u_n)_{n \in {\mathbb{Z}}}\,\in {\mathbb{C}}\,\,\,\vert \, || u||_{l^p}=\left(\sum_{n \in {\mathbb{Z}}}|u_n|^p\right)^{1/p}\right\}.
\end{equation}
When $p=2$, the Hilbert space $l^2$ is endowed with the inner product  
\begin{equation}
(u,v)_{l^2}=\sum_{n \in {\mathbb{Z}}}u_n \overline{v}_n,
\end{equation}
for any $u=(u_n)_{n \in {\mathbb{Z}}},v=(v_n)_{n \in {\mathbb{Z}}}\in l^2$  (where $\overline{v}_n$ denotes the conjugate of 
$v_n$) and induced norm
\begin{equation}
|| u||_{l^2}^2=(u,u).
\end{equation}
In some cases, when convenient, we shall use the inner product 
\begin{eqnarray*}
	(u,v)_{\ell^2}=\mathrm{Re}\sum_{{n\in\mathbb{Z}}}u_n\overline{v_n},\;\;u,\,v\in\ell^2,
\end{eqnarray*}
without a change in notation.  The continuous embeddings
\begin{equation}
l^r\subset l^s,\,\,\,|| w||_{l^s}\le || 
w||_{l^r},\,\,\,1 \le r\le s \le \infty,\label{eq:embeddings}
\end{equation}
will be thoroughly used in the sequel.

Next, let us recall some basic properties of the linear difference operators involved in the models. For any $u \in l^2$, the 1D-discrete Laplacian $\Delta_d:l^2\rightarrow l^2$,
\begin{equation}
(\Delta_d u)_{n }=u_{n+1}-2u_n+u_{n-1},
\end{equation}
is bounded, self-adjoint on $D(\Delta_d)=l^2$ and $\Delta_d\le 0$, due to the relations
\begin{equation}
\begin{split}
\label{prop}
(\Delta_{d} u,v)_{l^2}&=(u,\Delta_{d}v)_{l^2}, \quad \phi,\;\psi\in l^2,\\
(\Delta_{d}u,u)_{l^2}&=-\sum_{n\in\mathbb{Z}}|u_{n+1}-u_n|^2\leq 0,\\
||\Delta_d u||_{l^2}^2&\leq 4||u||_{l^2}^2.
\end{split}
\end{equation}
With the properties \eqref{prop}, it is established in \cite[Lemma 2.2, pg. 715]{NAE2007}, that the operator $\mathbf{L}:l^2\rightarrow l^2$, $(\mathbf{L}u)_{n\in\mathbb{Z}}=(1+i\alpha)(\Delta_d u)_{n\in\mathbb{Z}}$, generates a one parameter semigroup  $U(t)$ on $\ell^2$, as a perturbation of $-\Delta_d$ by the bounded (skew-adjoint) linear operator $A_1:\ell^2\rightarrow\ell^2$,  $(A_1u)_{n\in\mathbb{Z}}=-i\alpha(\Delta_d u)_{n\in\mathbb{Z}}$ \cite{Pazy83}. The semigroup $U(t)$  solves the underlying linear equation $\dot{u}=\mathbf{L}u$, namely $U(t)=\exp(\mathbf{L}t)$.
\paragraph{Local existence of solutions.}
With the preliminaries stated above, the initial value problem 
(\ref{eq:systemglobal})-(\ref{eq:icsglobal}) can be considered in the form of the integral equation
\begin{eqnarray}
\label{integ1}
u(t)=U(t)u^0+\int_{0}^{t}U(t-s)F(u(s))ds,\;\;\;u^0=\{u_{n,0}\}_{n\in\mathbb{Z}}.
\end{eqnarray}
%
%
where the operator $F$ is defined by
\begin{equation}
\label{nop}
\left(F(u)\right)_n=
(1-\delta)u_n+(1+i\alpha)(u_{n+1}+u_{n-1}-2u_n)-
(1+i\beta)\left[\gamma|u_n|^2u_n+\frac{\mu}{2}(u_{n+1}+u_{n-1})|u_n|^2\right]+g_n.
\end{equation}

In order to apply subsequently the local solvability theory to the integral equation \eqref{integ1},  the following lemma is a prerequisite:
%
\begin{lemma}
 \label{Lemma:Lipschitz} Let $g=(g_n)_{n\in {\mathbb{Z}}}\in l^2$ and $\alpha,\beta,\delta,\mu\in\mathbb{R}$. The operator $F$ defined \eqref{nop} is well defined as an operator $F:\,l^2\rightarrow l^2$ and 
is Lipschitz continuous on bounded sets of $l^2$.
\end{lemma}
\noindent{\bf Proof:} Let $u \in B_R$, with $B_R$ a closed ball of $l^2$ of center $0$ and radius $R$. We only need to consider
the nonlinear operator $N:l^2\rightarrow l^2$, $(N(u))_n=
-(1+i\beta)\left[\gamma|u_n|^2u_n+\frac{\mu}{2}(u_{n+1}+u_{n-1})|u_n|^2\right]$.  Recalling that $\sum_n|\theta_n|\,|\theta_{n+1}|\leq\sum_n|\theta_n|^2$, for any $\theta\in l^2$, we have:
\begin{eqnarray*}
 || N(u) ||_{l^2}^2&=&\sum_n \left|
 (1+i\beta)\left[\gamma|u_n|^2u_n+\frac{\mu}{2}(u_{n+1}+u_{n-1})|u_n|^2\right]\right|^2\nonumber\\
 &\leq &(1+\beta^2)\sum_n \left(|\gamma||u_n|^2|u_{n}|+\left|\frac{\mu}{2}\right||u_{n+1}+u_{n-1}||u_n|^2\right)^2
 \nonumber\\
 &=& \gamma^2(1+\beta^2)\sum_n|u_n|^4|u_n|^2+
 \frac{\mu^2}{4}(1+\beta^2)\sum_n |u_{n}|^4|u_{n+1}+u_{n-1}|^2+
 |\gamma| |\mu| \sum_n |u_{n}|^2|u_n||u_{n+1}+u_{n-1}||u_n|^2
  \nonumber\\
&\le& \gamma^2(1+\beta^2)\sup_{n}|u_{n}|^4\sum_n|u_n|^2+\mu^2(1+\beta^2)\sup_{n}|u_{n}|^4\sum_n|u_n|^2+2|\mu||\gamma|\sup_{n}|u_{n}|^4\sum_n|u_n|^2.
\end{eqnarray*}
Then, since $l^2\subset l^{\infty}$ with a continuous embedding, we derive the existence of a constant $K(\beta, \gamma, \mu, R)$, such that 
\begin{eqnarray*}
|| N(u) ||_{l^2}^2\leq	K(\beta, \gamma, \mu, R)||u||^2_{l^2}
\end{eqnarray*}
Hence, $N:\,l^2\rightarrow l^2$ is bounded on bounded sets of $l^2$. 
For $u,v \in B_R$ we have
\begin{eqnarray*}
\label{difN}
 || N(v)-N(u)||_{l^2}^2=
 \sum_{n \in {\mathbb{Z}}}\bigg|
 (1+i\beta)\left[\gamma|v_n|^2v_n+\frac{\mu}{2}(v_{n+1}+v_{n-1})|v_n|^2\right]-
 (1+i\beta)\left[\gamma|u_n|^2u_n+\frac{\mu}{2}(u_{n+1}+u_{n-1})|u_n|^2\right]
\bigg|^2.
\end{eqnarray*}
Note that in the difference $|u_n|^2u_n-|v_n|^2v_n$, when adding and subtracting simultaneously the term $|v_n|^2u_n$, we may write that 
\begin{eqnarray*}
\left||u_n|^2u_n-|v_n|^2v_n\right|\leq |v_n|^2|u_n-v_n|+|u_n|\left(|v_n|+|u_n|\right)|u_n-v_n|.
\end{eqnarray*}	
Similarly, for the difference $(v_{n+1}+v_{n-1})|v_n|^2-(u_{n+1}+u_{n-1})|u_n|^2$, when adding and subtracting simultaneously the term $(v_{n+1}+v_{n-1})|u_n|^2$, we may write
\begin{eqnarray*}
\left|(v_{n+1}+v_{n-1})|v_n|^2-(u_{n+1}+u_{n-1})|u_n|^2\right|&\leq& |v_{n+1}+v_{n-1}|\left(|v_n|+|u_n|\right)|v_n-u_n|\\
&&+|u_n|^2|v_{n+1}-u_{n+1}|+|u_n|^2|v_{n-1}-u_{n-1}|.	
\end{eqnarray*}	
Then, we may treat $|| N(v)-N(u)||_{l^2}^2$ with a similar algebra used above for the local boundedness of $N$ in $l^2$, to derive the existence of a constant $L(\beta,\gamma, \mu, R)$ such that 
\begin{eqnarray*}
|| N(v)-N(u)||_{l^2}^2\leq L(\beta,\gamma, \mu, R) ||v-u||_{l^2}^2,
\end{eqnarray*}	
verifying that the map $N:\,l^2\rightarrow l^2$ is Lipschitz continuous on bounded sets of $l^2$. \ \ $\square$

With Lemma \ref{Lemma:Lipschitz} in hand, the local existence result of solutions can be stated in the following proposition as an application of the local existence results for abstract evolution equations in Banach spaces \cite{Pazy83},\cite{Zeidler}. 
\begin{proposition} 
	\label{thloc}	
	Assume that $g_n\in l^2$, $\alpha,\beta, \delta, \mu\in\mathbb{R}$, and let $u^0\in l^2$, arbitrary. There exists some $T^*(u^0)>0$ such that the integral equation \eqref{integ1} has a unique solution $u\in C^1([0,T],l^2)$  for all $0<T<T^*(u^0)$. In addition, the following alternatives hold: Either $T^*(u^0)=\infty$ (global existence) or $T^*(u^0)<\infty$ and $\lim_{t\uparrow T^*(U^0)}||U(t)||_{\ell^2}=\infty$ (collapse). Furthermore the solution $u$ depends continuously on the initial condition $u^0\in l^2$, with respect to the norm of $C([0,T],l^2)$.
\end{proposition}
\paragraph{Global existence and asymptotic estimates.}
Proposition \ref{thloc} for the extended DGL (\ref{eq:systemglobal})-(\ref{eq:icsglobal}) establishes local existence of solutions for its limits,  the problems (\ref{eq:DNLS})-(\ref{eq:icsDNLS}) and (\ref{eq:ALS})-(\ref{eq:icsALS}).  For global existence of solutions we start with the case of the L-DGL (\ref{eq:DNLS})-(\ref{eq:icsDNLS}).
%
\begin{lemma}
\label{Lemma:unique}		
For the L-DGL (\ref{eq:DNLS}) assume $\delta>1$ and $\beta>0$. For every $u(0)=u^0\in l^2$, 
the system (\ref{eq:DNLS}) possesses a unique global solution $u(t)$ on $[0,\infty)$ belonging to
$C^1([0,\infty),l^2)$.
\end{lemma}
%
%
\noindent{\bf Proof:} 
For the L-DGL (1) we multiply the equation by $\overline{u}_n$, summing over $\mathbb{Z}$ keeping the real parts. The result is the energy balance equation:
\begin{eqnarray}
  \frac{d}{dt}\sum_{n\in {\mathbb{Z}}}|u_n|^2+2(\delta-1)||u||_{l^2}^2+2\sum_{n\in\mathbb{Z}}|u_{n+1}-u_n|^2=-2\sum_{n\in {\mathbb{Z}}} |u_n|^4+2\mathrm{Re}\sum_{n\in\mathbb{Z}}g_n\overline{u}_n.
  \label{eq:estimateattractorphi}
 \end{eqnarray}
  Using the Cauchy-Schwarz inequality and Young's inequality for the last term on the right-hand side of (\ref{eq:estimateattractorphi}) 
\begin{equation}
 \bigg|\sum_{n\in\mathbb{Z}}g_n\overline{u}_n\bigg|\leq \frac{(\delta-1)}{2}||u||^2_{l^2}+\frac{1}{2(\delta-1)}||g||^2_{l^2},
\end{equation}
 we obtain
\begin{equation}
  \frac{d}{dt}|| u(t)||_{l^2}^2 +(\delta -1)|| u(t)||_{l^2}^2\le \frac{1}{\delta-1}|| g||_{l^2}^2.\label{eq:DNLSl2above}
\end{equation}
By the assumption that $\delta>1$ it follows that $u \in L^{\infty}([0,\infty),l^2)$. 
With the use of Gronwall's inequality we obtain:
\begin{equation}
 || u(t)||_{l^2}^2\le || u(0)||_{l^2}^2\exp[-(\delta-1) t]+\frac{1}{(\delta-1)^2}|| g||_{l^2}^2\left(1-\exp[-(\delta-1) t]\right).
\end{equation}
Asymptotically for $t\rightarrow \infty$, yields the following upper bound:
\begin{equation}
 {\lim\sup}_{t \rightarrow \infty} || u(t)||_{l^2}^2\le \frac{1}{(\delta-1)^2}|| g||_{l^2}^2,
\end{equation}
implying the global existence of solutions for For the L-DGL (\ref{eq:DNLS})\ \ \ $\square$
\vspace*{0.5cm}

The parametric global existence regimes for the the NL-DGL (\ref{eq:ALS}) {\emph is much richer than for the L-DGL (\ref{eq:DNLS})}, as it is shown in the lemma below.
\begin{lemma}
	\label{nLDGL1}
	For the NL-DGL (\ref{eq:ALS}) we assume that $\delta>1$,  $\beta>0$, and consider the following cases for $\beta$, $||g||_{l^2}$ and the initial data $v(0)=v^0\in l^2$:
	\begin{enumerate}	
		\item
		Assume that $\delta$, $\beta$ and $||g||_{l^2}$  satisfy
		\begin{equation}
		\label{dALS0a}
		(\delta-1)^3<8\sqrt{1+\beta^2}||g||^2_{l^2}.
		\end{equation}
		Then, for any $v^0\in l^2$, the solutions exist globally in time and are uniformly bounded, satisfying
		\begin{equation}
		\label{eq:ALSbound}
		||v(t)||_{l^2}^2\leq ||v^0||_{l^2}^2.
		\end{equation}
		\item We consider the constants $B, R_1, R_2>0$ defined explicitly by the relations \eqref{dALS4}, \eqref{dALS5r} below. Assume that $\delta$, $\beta$ and $||g||_{l^2}$  satisfy
		\begin{equation}
		\label{dALS0b}
		(\delta-1)^3>8\sqrt{1+\beta^2}||g||^2_{l^2},
		\end{equation}
		and that the initial data $v^0\in l^2$ satisfy 
		\begin{eqnarray}
		\label{dALS0c}
		||v^0||_{l^2}^2\leq R_1.
		\end{eqnarray}
		Then, solutions exist globally for all $t>0$ and have the following asymptotic behaviour:
		\begin{eqnarray}
		\label{dANLS0d}
		R_2<||v(t)||_{l^2}^2\leq ||v^0||^2\leq R_1\;\;\mbox{for all}\;\; t > 0\;\;\mbox{and} \lim_{t\rightarrow\infty}||v(t)||_{l^2}^2=R_2.
		\end{eqnarray}
	\end{enumerate}
\end{lemma}
\noindent{\bf Proof:} We multiply equation \eqref{eq:ALS} by $\overline{v}_n$, sum over $\mathbb{Z}$, keeping the real parts, to get the following energy balance equation:
\begin{equation}
\label{dALS1}
\frac{1}{2}\frac{d}{dt}||v||_{l^2}^2+(\delta-1)||v||_{l^2}^2+\sum_{n\in\mathbb{Z}}|v_{n+1}-v_n|^2=-\frac{1}{2}\mathrm{Re}(1+i\beta)\sum_{n\in\mathbb{Z}}|v_n|^2(v_{n+1}+v_n)\overline{v}_n+\mathrm{Re}\sum_{n\in\mathbb{Z}} g_n\overline{v}_n.
\end{equation}
The right-hand side of \eqref{dALS1} can be estimated  from above, as follows:
\begin{equation}
\label{dALS2}
\bigg|\frac{1}{2}\mathrm{Re}(1+i\beta)\sum_{n\in\mathbb{Z}}|v_n|^2(v_{n+1}+v_n)\overline{v}_n\bigg|+\bigg|\sum_{n\in\mathbb{Z}}g_n\overline{v}_n\bigg|\leq \sqrt{1+\beta^2}||v||^4_{l^4}+\frac{(\delta-1)}{2}||v||_{l^2}^2+\frac{1}{2(\delta-1)}||g||^2_{l^2}.
\end{equation}
Then, combining \eqref{dALS1} and \eqref{dALS2}, we arrive at the differential inequality
\begin{equation}
\label{dALS3}
\frac{d}{dt}||v||_{l^2}^2+(\delta-1)||v||^2_{l^2}-2\sqrt{1+\beta^2}\big(||v||_{l^2}\big)^2-\frac{||g||^2_{l^2}}{\delta-1}\leq 0,
\end{equation}
which shall be explored further. For this purpose, we set $\chi(t)=||v(t)||^2_{l^2}$, and
\begin{equation}
\label{dALS4}
A=(\delta-1),\;\;B=2\sqrt{1+\beta^2},\;\;C=\frac{||g||^2_{l^2}}{\delta-1},
\end{equation}
and \eqref{dALS3} can be written as 
\begin{equation}
\label{dALS5}
\dot{x}(t)+A\chi(t)-B\chi(t)^2-C\leq 0,\;\;\;A,B,C>0.
\end{equation}

For further treatment of \eqref{dALS5} we distinguish between the following cases:
\begin{enumerate}
	\item $D=A^2-4BC<0$, which due to the definition of the constants $A,B,C$ in \eqref{dALS4} implies the condition \eqref{dALS0a}. In this case, the quantity $A\chi-B\chi^2-C>0$ for any $\chi\in\mathbb{R}$, and \eqref{dALS5} becomes $\dot{\chi}(t)\leq 0$, which gives the uniform bound \eqref{eq:ALSbound}.
	\item  $D=A^2-4BC>0$ which gives the condition \eqref{dALS0b}. Then the quadratic equation $A\chi-B\chi^2-C=0$ has two positive roots
	\begin{eqnarray}
	\label{dALS5r}
	R_{1,2}=\frac{(\delta-1)\pm\sqrt{(\delta-1)^2-\frac{8\sqrt{1+\beta^2}||g||^2_{l^2}}{(\delta-1)}}}{4\sqrt{1+\beta^2}}=\frac{A\pm\sqrt{D}}{2B}.
	\end{eqnarray} 
\end{enumerate}
Let now 
\begin{equation}
\label{dALS6}
\chi(t)=z(t)+R_{i},
\end{equation}
where $R_{i}$ is one of the roots $R_1$ or $R_2$. Since $R_i$ is a solution of the quadratic equation $A\chi-B\chi^2-C=0$, the $R_i$ satisfy inequality \eqref{dALS5} as an equality, that is
\begin{equation}
\label{dALS7}
\dot{R}_i+A R_i-B R_i^2-C=0,
\end{equation}
($\dot{R}_i=0$). Then, subtraction of \eqref{dALS7} from \eqref{dALS5},  yields that $z(t)$  satisfies the inequality \eqref{dALS5}. Thus, the change of variable defined by \eqref{dALS6} can be substituted into \eqref{dALS5}, and little algebra results in the Bernoulli inequality for $z(t)$:
\begin{equation}
\label{dALS9}
\dot{z}(t)+ (A-2BR_i)z(t)-Bz(t)^2\leq 0.
\end{equation}
To derive a uniform in time upper bound for the solutions, we need that the coefficient 
\begin{equation}
\label{dALS10}
A-2BR_i>0\Longrightarrow R_i<\frac{A}{2B}.
\end{equation}
Clearly, \eqref{dALS10} is satisfied only by the root $R_2$. Then, we set
\begin{eqnarray}
\label{dALS11}
K=A-2BR_2=\sqrt{D}>0.
\end{eqnarray}
Since \eqref{dALS11} holds for $K$, we treat the Bernoulli differential inequality as follows:  We use the change of variable $\psi=z^{-1}$. As $\dot{\psi}=-z^{-2}\dot z$, we have that $\dot{z}=-\dot{\psi}z^2$. Substituting in \eqref{dALS10}, we carefully integrating the latter as follows:
\begin{eqnarray}
\label{dALS12}
-\dot{\psi}z^2+Kz-Bz^2&\leq& 0\Rightarrow\nonumber\\
\dot{\psi}-K\psi+B&\geq& 0\Rightarrow\nonumber\\
\psi(t)&\geq & \bigg(\psi(0)-\frac{B}{K}\bigg)e^{Kt}+\frac{B}{K}.
\end{eqnarray}
We may ensure positivity of $\psi(t)$ for all $t\geq 0$, and thus of $z(t)$. This follows  by assuming that 
\begin{eqnarray*}
	\psi(0)\geq\frac{B}{K}\Longrightarrow z(0)\leq\frac{K}{B}\Longrightarrow \chi(0)\leq \frac{K}{B}+R_2=\frac{\sqrt{D}}{B}+R_2=R_1,
\end{eqnarray*}
giving the assumption on the initial data \eqref{dALS0c}. With this assumption on the initial data, we get from \eqref{dALS12} (when inverting for $z$), 
\begin{equation}
\label{dALS13}
z(t)\leq \bigg[\bigg(\psi(0)-\frac{B}{K}\bigg)e^{Kt}+\frac{B}{K}\bigg]^{-1}.
\end{equation}
Eventually from  \eqref{dALS13} we observe that 
\begin{eqnarray}
\label{dANLS14}
0<z(t),\;\;\mbox{for all}\;\;t>0,\;\;\mbox{and} \lim_{t\rightarrow\infty}z(t)=0.
\end{eqnarray} 
We use \eqref{dANLS14} to conclude with the description of the dynamics. Since $z(t)>0$ for all $t>0$, we have by its definition \eqref{dALS6}, that $\chi(t)=||v(t)||_{l^2}>R_2$. Thus, the solution belongs to the regime of positivity of the quadratic  $A\chi-B\chi^2-C$, for all $t>0$. With this  being the case, we derive again from \eqref{dALS5} that $\dot{x}(t)\leq 0$, for all $t>0$, implying the boundedness $||v(t)||^2_{l^2}\leq ||v(0)||^2_{l^2}$.  Finally, the definition of $z$ implies the limit for $||v(t)||_{l^2}^2$. \ \ $\square$\\
\\
In conclusion, for the L-DGL initial value problem (\ref{eq:DNLS})-(\ref{eq:icsDNLS}) solutions exist globally unconditionally with respect to the initial data, the forcing and the parameters of the system: for any  
$u^0 \in l^2$,  
the corresponding solution $u (t)$  of (\ref{eq:DNLS})-(\ref{eq:icsDNLS}) 
is bounded 
for all $t\in [0,\infty)$ and the  solution operator determined by 
\begin{equation}
 S_u(t):u^0 \in l^2 \rightarrow  u(t)=S_u(t)u^0\in l^2,\,\,\,t\ge 0,
\end{equation}
forms a continuous semigroup $\left\{ S_u(t)\right\}_{t\ge 0}$.
This is not the case for the solutions  $v(t)$ of the initial-value problem of the NL-DGL (\ref{eq:ALS})-(\ref{eq:icsALS}), as they exist globally for each of the in  Lemma \ref{nLDGL1}  classified regimes for the parameters, the forcing and the initial data.  Nevertheless, in these regimes the solution operator  of the NL-DGL problem (\ref{eq:ALS})-(\ref{eq:icsALS})
\begin{equation}
 S_v(t):v^0 \in l^2 \rightarrow  v(t)=S_v(t)v^0\in l^2,\,\,\,t\ge 0,
\end{equation}
generates a continuous semigroup $\left\{ S_v(t)\right\}_{t\ge 0}$ on $l^2$.
\section{Attractors and invariant sets}\label{section:attractors}
\setcounter{equation}{0}
In this section we establish the existence of restricted attractors, global attractors and other invariant sets for the dynamical systems associated to the discrete Ginzburg-Landau equations \eqref{eq:DNLS} and \eqref{eq:ALS}. 
\vspace{0.5cm}

\subsection{Existence of a restricted attractor for the semigroup   $\left\{S_v(t)\right\}_{t\ge 0}$ in $l^2$}
This section is devoted to the study of the semigroup $\left\{S_v(t)\right\}_{t\ge 0}$ of the NL-DGL system (\ref{eq:ALS})-(\ref{eq:icsALS}). The richer scenarios for invariant sets of the semigroup $\left\{S_v(t)\right\}_{t\ge 0}$ of the NL-DGL system  is already indicated by Lemma \ref{nLDGL1}. First, we discuss the implications of the  above possibilities on the domain of attraction and then we pay attention to the existence of a restricted attractor.

\subsubsection{Invariant sets in $l^2$ for the NL-DGL (\ref{eq:ALS})-(\ref{eq:icsALS})}
There exist at least two possibilities for attracting balls capturing the orbits of the semigroup  $\left\{S_v(t)\right\}_{t\ge 0}$ of the NL-DGL system (\ref{eq:ALS})-(\ref{eq:icsALS}), as well as the possibility of an invariant annulus. These scenarios correspond to distinct parametric regimes and have distinct domain of attractions as explained in the following lemma. 
%
%
\begin{figure}[tbp!]
	\begin{center}
		\begin{tabular}{cc}
			\includegraphics[width=.70\textwidth]{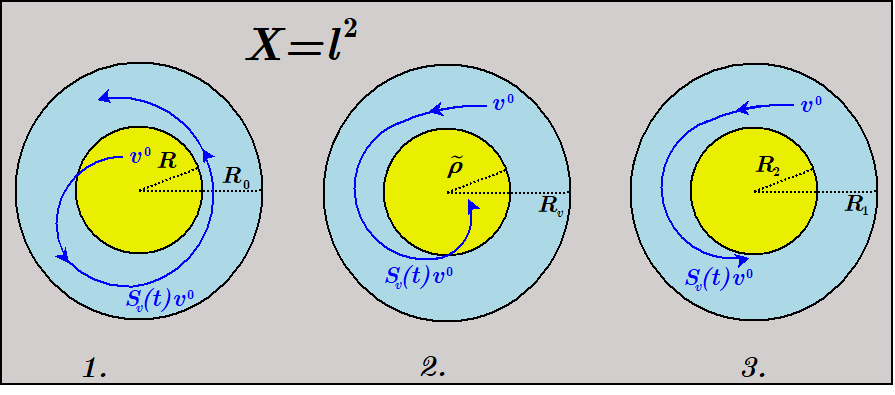}			
		\end{tabular}
	\end{center}
	\caption{The possibilities of Lemma \ref{Lemma:absorbingmu}. Case 1: The existence of a non-escaping ball $\mathcal{B}(0,R_0)$ with  $R_0(\delta, \beta, ||g||_{l^2}, R)$ for the orbit $S_v(t)v^0$ of the initial condition $v^0\in \overline{\mathcal{B}}(0,R)$. Case 2: The existence of a restricted absorbing ball $\mathcal{B}(0,\tilde{\rho})$ with a domain of attraction the closed ball $\overline{\mathcal{B}}(0,R_v)$. Case 3: The existence of the invariant annulus $\mathcal{A}_{1,2}$. For any initial condition $R_2<||v^0||_{l^2}^2\leq	R_1$, the orbit $v(t)=S_v(t)v^0$ satisfies $R_2<||v(t)||_{l^2}^2\leq	R_1$, and $lim_{t\rightarrow\infty}||v(t)||_{l^2}^2=R_2$. }
	\label{fig1}
\end{figure}
\begin{lemma}
\label{Lemma:absorbingmu}	
For the NL-DGL (\ref{eq:ALS}) we assume that $\delta>1$, $\beta>0$ and consider the following cases for $\beta$, $||g||_{l^2}$ and the initial data $v(0)=v^0\in l^2$:
	\begin{enumerate}	
		\item (Non-escaping ball).
		Assume that $\delta$, $\beta$ and $||g||_{l^2}$  are in the regime described by the condition \eqref{dALS0a} and that the initial data $v^0$ are in a closed ball $\overline{\mathcal{B}}(0,R)$ for some $R>0$.  
Then, there exists a closed ball $\overline{\mathcal{B}}(0,R_0)\subset l^2$ with $R_0(\delta, \beta, ||g||_{l^2}, R)$  for the semigroup $\left\{S_v(t)\right\}_{t\ge 0}$ of the system \eqref{eq:ALS}-\eqref{eq:icsALS},  with the property $\overline{\mathcal{B}}(0,R)\subset \overline{\mathcal{B}}(0,R_0)$, and
\begin{eqnarray}
\label{tas}
S_v(t) \overline{\mathcal{B}}(0,R)\subset \overline{\mathcal{B}}(0,R_0)\;\;\mbox{for all}\;\;t\geq 0.
\end{eqnarray}
\item (Absorbing ball for restricted initial data).	Assume that $\delta$, $\beta$ and $||g||_{l^2}$  are in the regime described by the condition \eqref{dALS0a} and that the initial data $v^0$ further satisfy
\begin{equation}
\label{rescond}
|| v^0||_{l^2}^2\le R_{v,0}^2:=\frac{\delta-1}{2\sqrt{1+\beta^2}},
\end{equation}
that is $v^0\in \mathcal{B}_{R_{v,0}}:=\overline{\mathcal{B}}(0,{R_{v,0}})\subset l^2$. 
Then, for the semigroup $\left\{S_v(t)\right\}_{t\ge 0}$
there exists a bounded absorbing set $\mathcal{B}_{v,0}$ in $l^2$ for the above restricted initial data, viz. for every bounded set $\mathcal{B}\subset \overline{\mathcal{B}}(0,{R_{v,0}})$ of $l^2$, there is a $t_0(\mathcal{B},\mathcal{B}_{v,0})$ s.t.  
$S_v(t)\mathcal{B}\subset\mathcal{B}_{v,0}\,\,\,\forall t\ge t_0(\mathcal{B},\mathcal{B}_{v,0})$.
\item (Positively invariant annulus in $l^2$).  Assume that $\delta$, $\beta$ and $||g||_{l^2}$ are in the regime described by the condition \eqref{dALS0b}. Then the annulus
\begin{eqnarray}
\label{inan}
\mathcal{A}_{1,2}=\bigg\{v\in l^2\;\;:\;\;R_2<||v||_{l^2}^2\leq R_1\bigg\},
\end{eqnarray}
is positively invariant for the semigroup $\left\{S_v(t)\right\}_{t\ge 0}$, viz.  for every bounded set $\mathcal{B}\subset \mathcal{A}_{1,2}$ it holds that $S_v(t)\mathcal{B}\subset \mathcal{A}_{1,2}$ for all $t\geq 0$. Moreover the limit  $\lim_{t\rightarrow\infty} S_v(t)v^0$ exists for any $v^0\in \mathcal{A}_{1,2}$, and   $\lim_{t\rightarrow\infty} S_v(t)v^0=R_2$.
\end{enumerate}
\end{lemma}


\noindent{\bf Proof:} 
\textit {1. (Non-escaping ball):} We turn again to the energy equation \eqref{dALS1}. By assumption, $||v^0||^2\leq R^2$. Then the right-hand side can be estimated as follows:
\begin{equation}
\label{dALSis1}
\begin{split}
\bigg| \mathrm{Re}(1+i\beta)|\sum_{n\in\mathbb{Z}}|v_n|^2(v_{n+1}+v_n)\overline{v}_n\bigg|&+2\bigg|\mathrm{Re}\sum_{n\in\mathbb{Z}}g_n\overline{v}_n\bigg|\\
&\leq 2\sqrt{1+\beta^2}||v||^4_{l^4}+(\delta-1)||v||^2_{l^2}+\frac{1}{\delta-1}|| g||_{l^2}^2\\
&\leq 2\sqrt{1+\beta^2}||v^0||^4_{l^2}+(\delta-1)||v||^2_{l^2}+\frac{1}{\delta-1}||g||_{l^2}^2.
\end{split}
\end{equation}
We set
\begin{equation}
\label{dALSis2}
\begin{split}
\rho_0&= 2\sqrt{1+\beta^2}||v^0||^4_{l^2}+\frac{1}{\delta-1}||g||_{l^2}^2,\\
\rho_0&\leq  2\sqrt{1+\beta^2}R^4+\frac{1}{\delta-1}||g||_{l^2}^2:=\rho_1(\beta, \delta,||g||_{l^2},R).
\end{split}
\end{equation}
Then, from \eqref{dALSis1} and \eqref{dALSis2}, the accordingly modified differential inequality obtained from equation \eqref{dALS1} implies the estimate
\begin{eqnarray}
\label{dALSis3}
|| v(t)||_{l^2}^2\le || v^0||_{l^2}^2\exp[-(\delta-1) t]+\frac{\rho_0}{\delta-1}|| g||_{l^2}^2\left(1-\exp[-(\delta-1) t]\right).
\end{eqnarray}
Letting $t\rightarrow \infty$ in \eqref{dALSis3}, we get
\begin{equation}
{\lim\sup}_{t \rightarrow \infty} || v(t)||_{l^2}^2\le \frac{\rho_0}{\delta-1}|| g||_{l^2}^2\leq \frac{\rho_1}{\delta-1}|| g||_{l^2}^2:=R_0^2,
\end{equation}
where $R_0$ is determined by the dependence of $\rho_1$ on the parameters $\beta$, $\delta$ and the initial radius $R$. Any orbit starting from the closed ball  $\overline{\mathcal{B}}(0,R)$ cannot escape from the closed ball $\overline{\mathcal{B}}(0,R_0)$ and each $\overline{\mathcal{B}}(0,R)$ defines a non-escaping ball $\overline{\mathcal{B}}(0,R_0)$. The latter is ``trivially attracting'' since \eqref{tas} holds, and the time of entry is $t_0=0$.\\
\textit {2. (Absorbing ball for restricted data):} The existence of the restricted absorbing ball follows by an alternative manipulation of the energy equation  \eqref{dALS1} and the estimate \eqref{dALSis1},  when we focus on the damping term on the right-hand side of the former. More precisely, instead of \eqref{dALSis1} we may have the following estimate: 
\begin{equation}
	\label{dALSis4}
	\begin{split}
		\bigg| \mathrm{Re}(1+i\beta)|\sum_{n\in\mathbb{Z}}|v_n|^2(v_{n+1}+v_n)\overline{v}_n\bigg|&+2\bigg|\mathrm{Re}\sum_{n\in\mathbb{Z}}g_n\overline{v}_n\bigg|\\
		&\leq 2\sqrt{1+\beta^2}||v||^2_{l^{\infty}}\sum_{n\in\mathbb{Z}}|v_n|^2+(\delta-1)||v||^2_{l^2}+\frac{1}{\delta-1}|| g||_{l^2}^2\\
		&\leq 2\sqrt{1+\beta^2}||v^0||^2_{l^2}\sum_{n\in\mathbb{Z}}|v_n|^2+(\delta-1)||v||^2_{l^2}+\frac{1}{\delta-1}||g||_{l^2}^2.
	\end{split}
\end{equation}
Now we assume that 
\begin{equation*}
	\label{rescon}
	\delta_0=(\delta -1)-2\sqrt{1+\beta^2}||v^0||^2_{l^2}>0,
\end{equation*}
which is actually condition \eqref{rescond}.
Inserting the modified estimate \eqref{dALSis4} in \eqref{dALS1} we get the differential inequality
\begin{equation*}
	\frac{d}{dt}||v||_{l^2}^2+\delta_0||v||_{l^2}^2\leq \frac{1}{\delta-1}||g||_{l^2}^2,
\end{equation*}
which implies the estimate 
\begin{equation*}
|| v(t)||_{l^2}^2\le || v(0)||_{l^2}^2\exp[-\delta_0 t]+\frac{1}{\delta_0(\delta-1)}|| g||_{l^2}^2\left(1-\exp[-\delta_0 t]\right).
\end{equation*}
Letting $t\rightarrow \infty$, yields that
\begin{equation}
\label{dALSis5}
{\lim\sup}_{t \rightarrow \infty} || v(t)||_{l^2}^2\le \frac{1}{\delta_0(\delta-1)}|| g||_{l^2}^2=\rho^2.
\end{equation}
Now, due to \eqref{dALSis5}, we observe  that for any number $\tilde{\rho}$, satisfying $\tilde{\rho}>\rho$, the ball $\mathcal{B}_{v,0}$ of $l^2$ centered at $0$ of radius $\tilde{\rho}$ is an absorbing set for the semigroup $S_v(t):$ That is, for a bounded set  
$\mathcal{B}\in l^2$  included in $\overline{\mathcal{B}}(0,{R_{v,0}})$, it follows that, for $t\ge t_0(\mathcal{B},\mathcal{B}_{v,0})$, where 
\begin{equation}
 t_0=\frac{1}{\tilde{\delta}}\log\left(\frac{R^2-\rho^2}{\tilde{\rho}^2-\rho^2}\right),
\end{equation}
 one has $|| v(t)||_{l^2}^2\le \tilde{\rho}^2$, i.e. 
$S_v(t)\mathcal{B}\subset \mathcal{B}_{v,0}$. \\
\textit{3. (Positively invariant annulus in $l^2$).} The proof is given in Lemma 	\ref{nLDGL1}, case 2.\ \ \ $\square$
\begin{remark}
	\label{remcon}
Regarding the scenarios of Lemma \ref{Lemma:absorbingmu}, note that scenario 1 may include scenario 2, since 1 holds under the most generic condition \eqref{eq:ALSbound} on the boundedness of the initial data and both scenarios hold under the same conditions on the parameters \eqref{dALS0a}. However, scenario 2 establishes under a further smallness condition on the initial data, namely \eqref{rescond}, the existence of the restricted absorbing ball $\mathcal{B}(0,\tilde{\rho})$ with an explicitly defined domain of attraction, the closed ball $\overline{\mathcal{B}}(0,R_v)$. On the other hand, scenario 3, is purely distinct since it holds when the parametric condition \eqref{dALS0a} is violated, that is, when \eqref{dALS0b} is satisfied and the initial data belong to the identified annulus $\mathcal{A}_{1,2}$.
\end{remark}
%
\subsubsection{Asymptotic compactness of the semigroup $\left\{S_v(t)\right\}_{t\ge 0}$}

Here, we verify that the semigroup $\left\{S_v(t)\right\}_{t\ge 0}$ associated with the NL-DGL (\ref{eq:ALS})-(\ref{eq:icsALS}) possesses the asymptotic tail-end property, if the hypotheses of Lemma \ref{Lemma:absorbingmu}.2 for the existence of the absorbing ball for the restricted initial data are fulfilled. We present an alternative approach for the proof instead of the tails-estimates method (see \cite{Bates, Zhou0,Zhou,Zhou1,Zhou3,Wang,Zhou2,Nikos,NAE2007,Du}). The approach is based on a contradiction argument and suitable implementation of the generalized comparison principle of ODE's \cite{McNabb1986}. We recall that the asymptotic tail property is essential in establishing the asymptotic compactness of the semigroup.

\begin{lemma}
\label{Lemma:asymtailmu}	
Assume the regime for the parameters, forcing and initial data stated in Lemma \ref{Lemma:absorbingmu}.2, that is, conditions (\ref{dALS0a}) and (\ref{rescond}). Let $(v_n(0))_{n\in {\mathbb{Z}}}=v^0 \in \mathcal{B}$ a bounded set of $l^2$, where $\mathcal{B}\subset \mathcal{B}_{R_v}$. Then, for any $\xi>0$, there exist
 $T(\xi)$ and $M(\xi)$ such that the solution $v(t)$ of (\ref{eq:ALS})-(\ref{eq:icsALS}) satisfies for all $t\ge T(\xi)$:
 \begin{equation}
  \sum_{|n|> 2K}|v_n(t)|^2\le \xi\,\,\,\,{\rm for\, any\,\,\,\,} K>M(\xi).\label{eq:asymptotic}
 \end{equation}
\end{lemma}


\noindent{\bf Proof:} For a contradiction, let us suppose that this assertion 
is not true, 
i.e. there is an $\epsilon_0>0$ and a subsequence $(n_k)_{k \in {\mathbb{Z}}}$ of $\mathbb{Z}$ s.t.
\begin{equation}
 \sum_{|k|>m}|v_{n_k}(t)|^2 \ge \epsilon_0,\,\,\,\forall t\ge 0,\,\,\,{\rm for\,\, any}\,\, m\in {\mathbb{N}}.\label{eq:epsilon0}
\end{equation}
Next, we elaborate the following argument based on the comparison principle \cite{McNabb1986}: Consider the equation
\begin{equation}
 \dot{w}(t)+Aw(t)-Bw^2(t)-C=0,\label{eq:upperw}
\end{equation}
with $A,B,C$ determined in (\ref{dALS4}) and $w(t)\ge 0$ for all $t\ge 0$.
Then, if $w(0)\ge\chi(0)$, by the comparison principle the solution $\chi (t)$ of (\ref{dALS5}) is bounded from above by the solution $w(t)$ of (\ref{eq:upperw}) for all 
$t \ge 0$.
Performing a Laplace transformation 
${\cal{L}}\{w(t)\}=\tilde{w}(s)$, yields
\begin{equation}
\label{lap1}
 s \tilde{w}(s)-w(0)+A \tilde{w}(s)=\frac{C}{s}+B F(s), 
\end{equation}
where $F(s)={\cal{L}}\{w^2(t)\}$. A rearrangement of the terms of \eqref{lap1}, gives the equation
\begin{equation}
 \tilde{w}(s)=\frac{1}{s+A}\left(w(0)+\frac{C}{s}+BF(s)\right).
\end{equation}
Applying the inverse Laplace transformation results in 
\begin{equation}
 w(t)=w(0)\exp(-At)+\frac{C}{A}\left(1-\exp(-At)\right)+B\int_0^t\exp(-A\theta)w^2(t-\theta)d\theta.
\end{equation}
Hence, we obtain the upper bound
\begin{eqnarray}
 \chi(t)&\le& \chi(0)\exp(-At)+\frac{C}{A}\left(1-\exp(-At)\right)+B\int_0^t\exp(-A\theta)\chi^2(t-\theta)ds\nonumber\\
 &\le& \chi(0)\exp(-At)+\frac{C}{A}\left(1-\exp(-At)\right)+B\max_{0\le \theta\le t}\chi^2(\theta)\int_0^t\exp(-A\theta)d\theta\nonumber\\
 &\le& \chi(0)\exp(-At)+\frac{1}{A}(C+B\chi^2(0))(1-\exp(-At),
\end{eqnarray}
where we used the relation (\ref{eq:ALSbound}). Asymptotically for $t\rightarrow \infty$, this gives the upper bound:
\begin{equation}
 \chi(t) \le \frac{1}{A}(C+B\chi^2(0)).\label{eq:upperchi}
\end{equation}

Now, let us note that since the NL-DGL equation is satisfied by the solution $v_n(t)$ for all $n\in\mathbb{Z}$, it is also satisfied by the subsequence $v_{n_k}(t)$, for all $n_k\in\mathbb{Z}$. We multiply equation \eqref{eq:ALS} by $\overline{v}_{n_k}$, sum over $|k|>m$, keeping the real parts to get the inequality
\begin{equation}
\label{baltailu}
  \frac{d}{dt}\sum_{|k|>m}|v_{n_k}(t)|^2 +(\delta-1) \sum_{|k|>m}|v_{n_k}(t)|^2 -2\sqrt{1+\beta^2}\left( \sum_{|k|>m}|v_{n_k}(t)|^2\right)^2-\frac{\sum_{|k|>m}|g_{n_k}|^2}{\delta-1}\leq 0,
\end{equation}
which is a differential inequality of the form (\ref{dALS5}) with constants 
\begin{equation*}
A=(\delta-1),\;\;B=2\sqrt{1+\beta^2},\;\;C=\frac{\sum_{|k|>m}|g_{n_k}|^2}{\delta-1}.
\end{equation*}
Since inequality (\ref{eq:upperchi}) is valid, one obtains
\begin{eqnarray}
 \sum_{|k|>m}|v_{n_k}(t)|^2 &\le& \sum_{|k|>m}|v_{n_k}(0)|^2 \exp(-(\delta-1)t)
 \nonumber\\
 &+&\frac{1}{\delta-1}\left(\frac{\sum_{|k|>m} |g_{n_k}|^2}{\delta-1}+2\sqrt{1+\beta^2}\left( \sum_{|k|>m}|v_{n_k}(0)|^2\right)^2
 \right)\left(1-\exp(-(\delta-1)t)\right).\label{eq:partialsum1}
 \end{eqnarray}
>From (\ref{eq:epsilon0}) and (\ref{eq:partialsum1}), we deduce that
 \begin{equation}
  \sup_{t\in[0,\infty]}\sum_{|k|>m}|v_{n_k}(t)|^2 =
  \frac{1}{(\delta-1)^2}\sum_{|k|>m} |g_{n_k}|^2 +2\frac{\sqrt{1+\beta^2}}{\delta-1}
  \left( \sum_{|k|>m}|v_{n_k}(0)|^2\right)^2
  \ge \epsilon_0,\,\,\,{\rm for\,\, any}\,\, m\in {\mathbb{N}}.
 \end{equation}
Since $g\in l^2$, by choosing $m$ sufficiently large,  this is equivalent to
\begin{equation}
 \sum_{|k|>m}|v_{n_k}(0)|^2
  \ge C(\beta,\delta)\tilde{\epsilon}>0,\label{eq:sumvnk}
  \end{equation}
  where 
 \begin{equation}
 \tilde{\epsilon}=\left(\epsilon_0  -
  \frac{1}{(\delta-1)^2}\sum_{|k|>m} |g_{n_k}|^2\right)^{1/2},\,\,\,\,\, C(\beta,\delta)= \left(\frac{\delta-1}{2\sqrt{1+\beta^2}}\right)^{1/2}.
 \end{equation}
 
 Thus, inequality (\ref{eq:sumvnk}) contradicts the fact that for $v^0 \in l^2$  every subsequence $(v_{n_k}(0))_{k\in {\mathbb{N}}}$ must  converge. \ \ \ $\square$
\\
\\
We also recall for the convenience of the reader the notion of the asymptotically compact semigroup.
\begin{definition} 
\label{def:asymcompact}	
The semigroup  $\left\{S_v(t)\right\}_{t\ge 0}$ is said to be asymptotically compact in $l^2$ if, for any bounded $\mathcal{B} \subset \mathcal{B}_{R_v}\subset  l^2$, and any sequence 
$\left\{t_n\right\}$,  $\left\{v_n\right\}$ with $t_n\ge 0$, $t_n \rightarrow \infty$ as $n\rightarrow \infty$, and $v_n \in \mathcal{B}$, the sequence  $\left\{S_v(t_n)v_n\right\} $ is relatively compact in $l^2$.
\end{definition}
With Lemma \ref{Lemma:asymtailmu} in hand, we establish the asymptotic compactness of the semigroup associated with the system (\ref{eq:ALS})-(\ref{eq:icsALS}) in the next proposition.
\begin{proposition}
{\it \,\,Under the same conditions of  Lemma \ref{Lemma:asymtailmu} the semigroup $\left\{S_v(t)\right\}_{t\ge 0}$  is asymptotically compact.}
\label{Proposition:asymcompmu}
\end{proposition}
\noindent{\bf Proof:} 
 By contradiction: Suppose that a subsequence of $v^n(t_n)=S_v(t_n)v^{0,n} \in \mathcal{B}$ diverges in $l^2$ as $t_n \rightarrow \infty$ (equivalently, the sequence $v^n(t_n)=S_v(t_n)v^{0,n}$ has not a convergent subsequence in $l^2$ as $t_n \rightarrow \infty$).
 Take any subsequence  $v^{n_k}(t_{n_k})=S_v(t_{n_k})v^{0,n_k} \in \mathcal{B}$, $k\in {\mathbb{N}}$. Suppose $||  v^{n_k}(t_{n_k}) ||_{l^2}^2 \rightarrow \infty$ as $t_{n_k}\rightarrow \infty$, viz. $k\rightarrow \infty$. 
 
 Then for any $M\in {\mathbb{R}}>0$ there are only finitely many $k$ s.t. 
 \begin{equation}
  ||  v^{n_k}(t_{n_k}) ||_{l^2}^2\le M.\label{eq:Mk}
 \end{equation}
 
  Denote all values of $k$ for which (\ref{eq:Mk}) is satisfied by $k_1,...,k_m$. Then, setting $N_M=\max\left\{k_1,...,k_m\right\}+1$, for any $k>N_M$ it holds that
  \begin{equation}
   ||  v^{n_k}(t_{n_k}) ||_{l^2}^2=\sum_{i\in {\mathbb{Z}}}|v^{n_k}_i(t_{n_{k}})|^2>M,\,\,\,\forall k>N_M.\label{eq:normpsiM}
  \end{equation}
 We  split the infinite sum in (\ref{eq:normpsiM}) as
 \begin{equation}
 \sum_{i\in {\mathbb{Z}}}|v^{n_k}_i(t_{n_{k}})|^2=\sum_{|i|>2L}|v^{n_k}_i(t_{n_{k}})|^2+\sum_{|i|\le 2L}|v^{n_k}_i(t_{n_{k}})|^2
 >M,\,\,\,\forall k>N_M,\label{eq:normpsiM1}
  \end{equation}
  for any fixed $0<L<\infty$. For the finite sum, we get $\sum_{|i|\le 2L}|v^{n_k}_i(t_{n_{k}})|^2<M_L$ with $M>M_L$, implying that
  \begin{equation}
 \sum_{|i|\ge 2L}|v^{n_k}_i(t_{n_{k}})|^2
 >M-M_L>0,\,\,\,\forall k>N_M.\label{eq:normpsiM2}
  \end{equation}
  Since the relation (\ref{eq:normpsiM2}) holds for all 
  $t\ge t_{n_{k> N_M}}$, for all $M>0$ and any $L>0$, it contradicts the asymptotic tail end property of $\left\{S_v(t)\right\}_{t\ge 0}$  as established by Lemma \ref{Lemma:asymtailmu}. 
  
  \ \ \ $\square$

\vspace*{0.5cm}

Finally, facilitating Proposition \ref{Proposition:asymcompmu} and Theorem 1.1 in \cite{Temam}, we are now able to state the main result of this section.

\begin{theorem}
 {\it \,\,
 The semigroup $\left\{S_v(t)\right\}_{t\ge 0}$ associated with the NL-DGL  (\ref{eq:ALS})-(\ref{eq:icsALS}), possesses a unique restricted global attractor ${\cal A}_v\subset \mathcal{B}_{v,0}\subset l^2$.}\label{Theorem:attractorALS}
 \end{theorem}
 
\subsection{Existence of a global attractor for the semigroup $\left\{S_u(t)\right\}_{t\ge 0}$}\label{section:Au}
Here, we prove the existence of an absorbing set in $l^2$ for the dynamical system belonging to the L-DGL (\ref{eq:DNLS})-(\ref{eq:icsDNLS}) in the asymptotic regime $t \rightarrow \infty$.

 
\begin{lemma}
Assume that  $\delta>1$, $\beta>0$ and let $(g_n)_{n\in {\mathbb{Z}}}=g \in l^2$ and $(u_n(0))_{n\in {\mathbb{Z}}}=u^0$. For the dynamical system determined by (\ref{eq:DNLS})-(\ref{eq:icsDNLS})
\begin{equation}
 S_u(t):\,u^0\in l^2\rightarrow u(t)\in l^2
\end{equation}
there exists  a bounded absorbing set $\mathcal{B}_{u,0}$ in $l^2$, viz. for every bounded set $\mathcal{B}$ of $l^2$, there is a $t_0(\mathcal{B},\mathcal{B}_{u,0})$ s.t.  
$S_u(t)\mathcal{B}\subset \mathcal{B}_{u,0}\,\,\,\forall t\ge t_0(\mathcal{B},\mathcal{B}_{u,0})$. 
\label{Lemma:asymboundgamma}
\end{lemma}


\noindent{\bf Proof:} 
Defining 
\begin{equation}
 \rho^2=\frac{1}{(\delta-1)^2}|| g||_{l^2}^2,
\end{equation}
we observe from Lemma \ref{Lemma:unique}, that for any number $\tilde{\rho}>\rho$ the ball $\mathcal{B}  _{u,0}\in l^2$  centered at $0$ of radius $\tilde{\rho}$ is an absorbing set for the semigroup $S_u(t).$ That is, if $\mathcal{B}$ is a bounded set of $l^2$ included in a ball centered at $0$ of radius $R$, it follows that for $t\ge t_0(\mathcal{B},\mathcal{B}_{u,0})$ where 
\begin{equation}
 t_0=\frac{1}{\delta-1}\log\left(\frac{R^2-\rho^2}{\tilde{\rho}^2-\rho^2}\right),
\end{equation}
one has $|| u(t)||_{l^2}^2\le \tilde{\rho}^2$, viz.
$S_u(t)\mathcal{B}\subset \mathcal{B}_{u,0}$.  \ \ \ $\square$


\vspace*{0.5cm}

Concerning the asymptotic tail-end property, for 
the L-DGL (\ref{eq:DNLS})-(\ref{eq:icsDNLS}) we have the following Lemma. 

\vspace*{0.5cm}
\begin{lemma}
 {\it \,\,  Let $(u_n(0))_{n\in {\mathbb{Z}}}=u^0 \in \mathcal{B}$, where $\mathcal{B}$ is a bounded set of $l^2$ and 
 $(g_n)_{n\in {\mathbb{Z}}}=g \in l^2$. For any $\xi>0$ there exist
 $T(\xi)$ and $M(\xi)$ s.t. the solution $u(t)$ of (\ref{eq:DNLS}),(\ref{eq:icsDNLS}) satisfies for all $t\ge T(\xi)$:
 \begin{equation}
  \sum_{|n|\ge 2K}|u_n(t)|^2\le \xi\,\,\,\,{\rm for\, any\,\,\,\,} K>M(\xi).\label{eq:asymptoticDNLS}
 \end{equation}}\label{Lemma:asymtailgamma}
 \end{lemma}

 \noindent{\bf Proof:} For a contradiction, let us suppose that this assertion is 
is not true, 
i.e. there is an $\epsilon_0>0$ and a subsequence $(n_k)_{k \in {\mathbb{Z}}}$ of $\mathbb{Z}$ s.t.
\begin{equation}
 \sum_{|k|>m}|u_{n_k}(t)|^2 \ge \epsilon_0,\,\,\,\forall t\ge 0,\,\,\,{\rm for\,\, any}\,\, m\in {\mathbb{N}}.
\end{equation}
Working as for the derivation of \eqref{baltailu}, we get the differential inequality for $u_{n_k}$,
\begin{equation}
  \frac{d}{dt}\sum_{|k|>m}|u_{n_k}(t)|^2 +(\delta-1) \sum_{|k|>m}|u_{n_k}(t)|^2\le \frac{1}{\delta-1}  \sum_{|k|>m}|g_{n_k}|^2,
\end{equation}
and Gronwall's inequality gives
\begin{equation}
 \sum_{|k|>m}|u_{n_k}(t)|^2 \le \sum_{|k|>m}|u_{n_k}(0)|^2 \exp(-(\delta-1)t)+\frac{\sum_{|k|>m} |g_{n_k}|^2}{(\delta-1)^2}\left(1-\exp(-(\delta-1)t)\right).\label{eq:partialsum}
 \end{equation}
 Conclusively,
 \begin{equation}
  \inf_{t\in[0,\infty]}\sum_{|k|>m}|u_{n_k}(t)|^2 =
  \frac{\sum_{|k|>m} |g_{n_k}|^2}{(\delta-1)^2} \ge \epsilon_0>0,\,\,\,{\rm for\,\, any}\,\, m\in {\mathbb{N}}.
 \end{equation}
 Therefore, for every $m\in {\mathbb{N}}$, $\sum_{|k|>m}|g_{n_k}|^2>C(\delta)\epsilon_0>0$,  which contradicts the fact that for $(g_n)_{n\in {\mathbb{Z}}}=g$ with $|| g  ||_{l^2}<\infty$, every subsequence $(g_{n_k})_{k\in {\mathbb{N}}}$ must  converge. 
 
 \ \ \ $\square$
\\
\\
 



\begin{proposition}
{\it \,\,Under the same conditions of Lemma \ref{Lemma:asymtailgamma},  the semigroup $\left\{S_u(t)\right\}_{t\ge 0}$ is asymptotically compact.}
\label{Proposition:asymcompgamma}
\end{proposition}

The proof of 
Proposition \ref{Proposition:asymcompgamma}  proceeds in a similar manner to the corresponding proof for the NL-DGL and is omitted here.

\vspace*{0.5cm}

In conclusion, by virtue of  Proposition \ref{Proposition:asymcompmu} and Theorem 1.1 in \cite{Temam}, we have:

\begin{theorem}
{\it \,\,
The semigroup $\left\{S_u(t)\right\}_{t\ge 0}$ attributed to the L-DGL  (\ref{eq:DNLS})-(\ref{eq:icsDNLS}) has  a unique global attractor ${\cal A}_u\subset \mathcal{B}_{u,0}\subset l^2$.}\label{Theorem:attractorDNLS1}
\end{theorem}
\section{Closeness of the NL-DGL and L-DGL solutions and congruence of the attractors}\label{section:closeness}
\setcounter{equation}{0}
\subsection{Closeness of the solutions}
In this section, we give the proof of Theorem \ref{Theorem:closeness}. The result gives in a sense of  ``continuous dependence'' the closeness of the solutions of  (\ref{eq:DNLS}) and  (\ref{eq:ALS}) for small enough initial data. A first application can be found in establishing the persistence of small amplitude dissipative solitons for the L-DGL \eqref{eq:DNLS}, close to those identified for the NL-DGL \eqref{eq:ALS}, when the suitable parametric regimes for the models described in \cite{AGL3,AGL4,UGL6} are considered.
\\
\\
\noindent{\bf Proof of Theorem \ref{Theorem:closeness}:} \textit{1. Case $\delta>1$}. One derives for $t>0$, the following differential inequality:
\begin{eqnarray}
\label{dif1}
  \frac{d}{dt}|| \Delta||_{l^2}^2&=&2|| \Delta||_{l^2} \frac{d}{dt}|| \Delta||_{l^2}\nonumber\\
&=&  \sum_{n \in {\mathbb{Z}}}\bigg\{ 2(1-\delta) |\Delta_n|^2+[(1+i\alpha)(\Delta_{n+1}+\Delta_{n-1}-2\Delta_n)\overline{\Delta}_n  +c.c.]\nonumber\\
&-&[(1+i\beta)|u_n|^2u_n \overline{\Delta}_n  +c.c.]\nonumber\\
&-&\left[\frac{1}{2}(1+i\beta)|v_n|^2(v_{n+1}+v_{n-1})\overline{\Delta}_n   
 +c.c.\right]\bigg\}\nonumber\\
&\le& 2(1-\delta) || \Delta||_{l^2}^2+
2(1+\beta)( \sup_{n \in {\mathbb{Z}}}|u_n|^2 || u||_{l^2} + 
\sup_{n \in {\mathbb{Z}}}|v_n|^2 || v||_{l^2})|| \Delta||_{l^2} \nonumber\\
   &\le& 2(1-\delta) || \Delta||_{l^2}^2+
2(1+\beta)( || u||_{l^2}^3 + 
 || v||_{l^2}^3)|| \Delta||_{l^2} 
\end{eqnarray}
For the derivation of \eqref{dif1}, we exploited the continuous embeddings \eqref{eq:embeddings}.
For $t>0$, it holds that 
\begin{equation}
\label{dif2}
 \frac{d}{dt}|| \Delta||_{l^2}+ 2(\delta-1) || \Delta||_{l^2}\leq 2(1+\beta)(R_u^3+R_v^3).
\end{equation}
Applying Gronwall's inequality to \eqref{dif2}, one gets
\begin{equation}
 \label{dif3}
 ||\Delta(t)||_{l^2}\leq ||\Delta (0)||\exp[-2(\delta-1)t]+\frac{(1+\beta)}{\delta-1}(R^3_u+R^3_v)\big[1-\exp[-2(\delta-1)t]\big].
\end{equation}
From \eqref{dif3}, we see that the estimate \eqref{eq:boundy} is satisfied with the 
constant 
\begin{equation}
\label{cdc0}
C=C_0+\bigg[\frac{(1+\beta)}{\delta-1}(C_{u,0}^3+C_{v,0}^3)\bigg].
\end{equation}
Taking the limit as $t\rightarrow\infty$,  in \eqref{dif3}, we deduce that
\begin{eqnarray*}
\limsup_{t\rightarrow\infty}||\Delta(t)||_{l^2}\leq \bigg[\frac{(1+\beta)}{\delta-1}(C_{u,0}^3+C_{v,0}^3)\bigg]\epsilon^3.
\end{eqnarray*}
justifying \eqref{cd1} with
\begin{eqnarray}
\label{cdc}
C_1=\bigg[\frac{(1+\beta)}{\delta-1}(C_{u,0}^3+C_{v,0}^3)\bigg]\epsilon^3.
\end{eqnarray}
\textit{2. Case $\delta\leq 1$.}  When $\delta=1$, we derive the claimed estimate from \eqref{dif2}, with the constant 
\begin{equation*}
C_2=C_0+2(1+\beta)(C_{u,0}^3+C_{v,0}^3)T_f.
\end{equation*}
Furthermore, when $\delta<1$, we get from \eqref{dif2} the estimate
\begin{equation}
	\label{dif40}
\begin{split}
||\Delta(t)||_{l^2}&\leq ||\Delta (0)||\exp[2(1-\delta)t]+\frac{(1+\beta)}{1-\delta}(R^3_u+R^3_v)\big[\exp[-2(1-\delta)t]-1\big],
\end{split}
\end{equation}
which holds for any $t\in (0,\infty]$.  Hence, for any $0<T_f<\infty$, by setting  this time
\begin{equation}
\label{dif4}
C_2=C_0\exp[2(1-\delta)T_f]+\frac{(1+\beta)}{1-\delta}(C_{u,0}^3+C_{v,0}^3)\big[\exp[2(1-\delta)T_f]-1\big],
\end{equation}
we conclude with the claimed  estimate \eqref{eq:boundylin}. \ \  $\square$

After having shown that the distance between the solutions of the 
NL-DGL and the L-DGL measured in terms of the $l^2-$metric remains small (bounded above by ${\cal{O}}(\epsilon^3)$), compared to the $l^2-$norm of the solutions themselves (bounded above by ${\cal{O}}(\epsilon$), an analogous feature directly follows for the $l^{\infty}-$norm (sup norm) 
determining the maximal distance between individual units, which is a consequence of the embeddings \eqref{eq:embeddings}. We have  the following
\begin{corollary}
	\label{infdist}
Assume the conditions of Theorem \ref{Theorem:closeness}. Then the closeness estimates \eqref{eq:boundy}-\eqref{cd1} for $\delta>1$, and \eqref{eq:boundylin} for $\delta\leq 1$, are satisfied by the $l^{\infty}$-distance of solutions $||\Delta(t)||_{l^{\infty}}=||u(t)-v(t)||_{l^{\infty}}$.	
\end{corollary}
\subsection{Congruence of the attractors ${\cal A}_v$ and ${\cal A}_u$}\label{section:congruence}
A second application of Theorem \ref{Theorem:closeness} is that it can be used in the dissipative regime $\delta>1$ to establish the congruence of the attractors ${\cal A}_u$ and ${\cal A}_v$ associated with the dynamical systems $S_u(t)$ and $S_v(t)$, that is the proof of Theorem \ref{Theorem:congruence}. \\
\\
\noindent{\bf Proof of Theorem \ref{Theorem:congruence}:}
For any bounded subset $\mathcal{B}\subset l^2$, it holds that 
\begin{eqnarray}
 {\rm dist}({\cal A}_v,{\cal A}_u)&\le& {\rm dist}({\cal A}_v,S_v(t)\mathcal{B})+   {\rm dist}(S_v(t)\mathcal{B},S_u(t)\mathcal{B})+ {\rm dist}(S_u(t)\mathcal{B},{\cal A}_u).
\end{eqnarray}
  As ${\cal A}_v$ attracts any bounded set $\mathcal{B}\subseteq \mathcal{B}_{R_v}\subset l^2$, for  
    any $\xi >0$, there is some $T_v(\xi)>0$ s.t. 
  \begin{equation}
  {\rm dist}\left({\cal A}_v,S_v(t)\mathcal{B}\right)=\sup_{a\in 
 {\cal A}_v}\,\inf_{v^0 \in \mathcal{B}} {\rm dist}\left(a,S_v(t)v^0\right)_{l^2} <\frac{\xi}{3},\,\,\,\forall t\ge T_v. \label{eq:ineqc1}
 \end{equation}
  Analogously, as   ${\cal A}_u$ attracts any bounded set $\mathcal{B}\subset l^2$ (${\cal A}_u$ actually attracts any bounded set in $l^2$ anyway),  
  for 
    any $\xi >0$, there is some $T_u(\xi)>0$ s.t. 
  \begin{equation}
  {\rm dist}\left(S_u(t)\mathcal{B},{\cal A}_u\right)={\rm dist}\left({\cal A}_u,S_u(t)\mathcal{B}\right)
  =\sup_{a\in 
 {\cal A}_u}\,\inf_{u^0 \in \mathcal{B}} {\rm dist}\left(a,S_u(t)u^0\right)_{l^2}
  <\frac{\xi}{3},\,\,\,\forall t\ge T_u. \label{eq:ineqc2}
 \end{equation}
 Let $\overline{T}=\max\{T_v(\xi),T_u(\xi)\}$ and consider the time interval $[0,t_0]$ with $t_0\ge \overline{T}$.
In light of Theorem \ref{Theorem:closeness}-case 1, for $\delta>1$, we have that for  every $t_0>\overline{T}$, there exists a constant $C>0$ s.t. for all 
$v^0,\,u^0 \in \mathcal{B} \subseteq \mathcal{B}_{R_v}$, fulfilling hypotheses \eqref{eq:initialuv}, it holds that for every $t\in [0,t_0]$,
 \begin{equation}
{\rm dist}(S_v(t)\mathcal{B},S_u(t)\mathcal{B})=\inf_{v^0 \in \mathcal{B}}  {\rm dist}_{l^2}\left(S_v(t)v^0,S_u(t)v^0\right)
  <C{\epsilon^3},\,\,\,\forall \epsilon>0.\label{eq:ineqc3}
 \end{equation}

 Combining (\ref{eq:ineqc1}),(\ref{eq:ineqc2}) and (\ref{eq:ineqc3})
we get 
\begin{eqnarray}
  0\le {\rm dist}\left({\cal A}_v,{\cal A}_u\right)&\le& 
  \frac{\xi}{3}+C\epsilon^3+\frac{\xi}{3}=\frac{2\xi}{3}+C\epsilon^3,\,\,\,\forall \xi>0.
\end{eqnarray}
By setting $\epsilon_0^3=\xi/(3C)$, one gets 
\begin{equation}
 {\rm dist}\left({\cal A}_v,{\cal A}_u\right)< \xi,\,\,\,\forall \epsilon\in [0,\epsilon_0],  
\end{equation}
from which, by the arbitrariness of $\xi$ it follows that 
${\rm dist}\left({\cal A}_v,{\cal A}_u\right)=0$,
and the proof is finished.\ \ $\square$

\section{Conclusions}
\label{SecCon}
While the integrable, non-local Ablowitz-Ladik (AL) lattice \cite{Ablowitz} is less significant in physical applications than its local variant-the  celebrated Discrete Nonlinear Schr\"odinger equation (DNLS), its counterpart incorporating linear and nonlinear damping (or even forcing), namely the non-local Discrete Ginzburg-Landau equation, has proven to be more relevant in various physical contexts particularly in nonlinear optics \cite{AchBook} and  in the mathematical description of arrays of vortices in fluids \cite{AGL1}-\cite{AGL2}.  For this important model, we have identified distinct scenarios for the asymptotic behavior of its globally existing solutions, depending on the involved parameters and external driving. 
Crucially, we have proven closeness of the  solutions of the non-local Discrete Ginzburg-Landau equation to the ones of its local counterpart; the latter  is another physically relevant nonlinear lattice. As a consequence, the compact attractors governing the asymptotic dynamics coincide when the initial data are selected within a suitable domain of attraction defined by the restricted absorbing ball of the non-local system.  
Our results on the  closeness of globally existing solutions and the congruence of attractors may inspire further studies in the context of dissipative localized structures and their dynamics taking into account stochastic and other random effects \cite{Stoch1}-\cite{Stoch4} or  space and time discretisations \cite{GLord},\cite{Han} relevant for the  numerical analysis of the (continuum) Ginzburg-Landau partial differential equation. 
One direction of research shall be the investigation of the inviscid limit \cite{Zhao} of the dissipative lattices (such as the discrete GLEs considered in this paper) to the Hamiltonian ones (represented by the integrable Ablowitz-Ladik lattice with its variety of analytical solutions)  motivated also by current closeness results for the conservative AL and DNLS lattices \cite{DNJ2021}. Corresponding investigations are in progress and will be reported in future publications. 
\section*{Acknowledgment}
We would like to thank the referee for his/her constructive comments and suggestions. 

\end{document}